# ASYMPTOTIC ENTROPY AND GREEN SPEED FOR RANDOM WALKS ON COUNTABLE GROUPS

By Sébastien Blachère , Peter Haïssinsky
and Pierre Mathieu

*Eurandom, Université Aix-Marseille 1
and Université Aix-Marseille 1*

We study asymptotic properties of the Green metric associated with transient random walks on countable groups. We prove that the rate of escape of the random walk computed in the Green metric equals its asymptotic entropy. The proof relies on integral representations of both quantities with the extended Martin kernel. In the case of finitely generated groups, where this result is known (Benjamini and Peres [*Probab. Theory Related Fields* **98** (1994) 91–112]), we give an alternative proof relying on a version of the so-called fundamental inequality (relating the rate of escape, the entropy and the logarithmic volume growth) extended to random walks with unbounded support.

**1. Introduction.** Let $\Gamma$ be an infinite countable group and let $(Z_n)$ be a transient random walk on $\Gamma$.

In order to study asymptotic properties of the random walk, we define the Green (or hitting) metric,

$$d_G(x,y) = -\ln \mathbb{P}^x[\tau_y < \infty],$$

where $\tau_y$ is the hitting time of the element $y$ by the random walk started at $x$.

Looking at the random walk via the Green metric leads to a nice geometrical interpretation of probabilistic quantities describing the long-time behavior of the walk. We illustrate this claim by showing that the rate of escape computed in $d_G$ coincides with the asymptotic entropy of the random walk; see Theorem 1.1. As another example of interest of the Green metric, we also explain how the Martin compactification of $\Gamma$ can be interpreted as the Busemann compactification of $\Gamma$ equipped with $d_G$. In a forthcoming paper [5], we use the Green metric to study fine geometric properties of the harmonic measure on boundaries of hyperbolic groups.

Before stating our theorem, let us first recall some definitions. The rate of escape of the random walk computed in the Green metric (in short, the *Green speed*) is defined by the almost sure limit

$$\ell_G \stackrel{\text{def.}}{=} \lim_{n \to \infty} \frac{d_G(e, Z_n)}{n}.$$

The asymptotic entropy of the random walk is defined by

$$h \stackrel{\text{def.}}{=} \lim_{n \to \infty} \frac{-\ln \mu^n(Z_n)}{n},$$

where $\mu$ is the law of the increment of the random walk (i.e., the law of $Z_1$) and $\mu^n$ is the $n$th convolution power of $\mu$ (i.e., the law of $Z_n$). This limit exists almost surely and is finite if the entropy of $\mu$,

$$H(\mu) \stackrel{\text{def.}}{=} -\sum_{x \in \Gamma} \mu(x) \ln \mu(x),$$

is finite. The asymptotic entropy $h$ plays a very important role in the description of the long-time behavior of the random walk, as illustrated in Derriennic [9, 10], Guivarc'h [15], Kaimanovich [17], Kaimanovich and Vershik [18] and Vershik [24], among others. For instance, it is known that $h = 0$ if and only if the Poisson boundary of the random walk is trivial.

Our main result is the following.

THEOREM 1.1. *For any transient random walk on a countable group such that $H(\mu) < \infty$, the asymptotic entropy $h$ and the Green speed $\ell_G$ are equal.*

In Section 2, we prove this result using an integral representation of $h$ on the Martin boundary of $\Gamma$ (Lemma 2.6) and interpreting the Green speed of the random walk as a limit of a Martin kernel (Proposition 2.4). This proof does not use any quantitative bound on the transition probabilities of the random walk and therefore applies to transient random walks on all countable groups, even those which are not finitely generated.

In Section 3, we consider the case of a finitely generated group $\Gamma$ and discuss the connection of Theorem 1.1 with the so-called "fundamental inequality" $h \leq \ell \cdot v$, where $\ell$ and $v$ denote the rate of escape and the logarithmic volume growth in some left-invariant metric on the group with a finite first moment. We first derive a new general version of the fundamental inequality for any random walk (with bounded or unbounded support) and any (geodesic or nongeodesic) left-invariant metric on the group with a finite first moment; see Proposition 3.4. We then use heat kernel estimates to obtain bounds on the logarithmic volume growth in the Green metric; see Proposition 3.1. Thus, we finally obtain another proof of Theorem 1.1, valid for finitely generated groups of superpolynomial volume growth. In the case of groups with polynomial volume growth, $h$ and $\ell_G$ are both zero.

For finitely generated groups, Benjamini and Peres [3] gave a different proof of the equality $h = \ell$. Even if their proof is written for finitely supported random walks, their method also works for random walks with infinite support (see the proof of Proposition 3.1).

## 2. Countable groups.

2.1. *The Green metric.* We will give the definition of the Green metric associated with a transient random walks and recall some of its properties from Blachère and Brofferio [4].

Let $\mu$ be a probability measure on $\Gamma$ whose support generates the whole group $\Gamma$. (We will always make this generating hypothesis.) We do not assume that $\mu$ is symmetric nor that it is finitely supported. Let $(X_k)$ be a sequence of i.i.d. random variables whose common law is $\mu$. The process

$$Z_k \stackrel{\text{def.}}{=} xX_1X_2\cdots X_k,$$

with $Z_0 = x \in \Gamma$, is an irreducible random walk on $\Gamma$ starting at $x$ with law $\mu$. We denote by $\mathbb{P}^x$ and $\mathbb{E}^x$, respectively, the probability and expectation related to a random walk starting at $x$. When $x = e$ (the identity of the group), the exponent will be omitted.

From now on, we will always assume the random walk to be transient, that is, with positive probability, it never returns to its starting point. This assumption is always satisfied if $\Gamma$ is not a finite extension of $\mathbb{Z}$ or $\mathbb{Z}^2$ (see Woess [25], Section I.3.B). On a finite extension of $\mathbb{Z}$ or $\mathbb{Z}^2$, there exists a canonical projection $\varphi$ onto an Abelian subgroup ($\{e\}, \mathbb{Z}$ or $\mathbb{Z}^2$); see Alexopoulos [1]. We define the first moment of the canonical projection of the random walk,

$$M_1(\mu) \stackrel{\text{def.}}{=} \sum_{x \in \Gamma} \|\varphi(x)\|\mu(x),$$

where $\|\varphi(x)\|$ is the norm of $\varphi(x)$. When $M_1(\mu) < \infty$, the random walk is transient if and only if it has a nonzero drift $[\sum_{x \in \Gamma} \varphi(x)\mu(x) \neq 0]$. But there are examples of recurrent and transient random walks with $M_1(\mu) = \infty$. There are even examples of transient symmetric random walks on $\mathbb{Z}$. For these results and examples, see Spitzer [23].

The *Green function* $G(x,y)$ is defined as the expected number of visits at $y$ for a random walk starting at $x$:

$$G(x,y) \stackrel{\text{def.}}{=} \mathbb{E}^x\left[\sum_{n=0}^{\infty} \mathbb{1}_{\{Z_n=y\}}\right] = \sum_{n=0}^{\infty} \mathbb{P}^x[Z_n = y].$$

Since the random walk is chosen to be transient, the Green function is finite for every $x$ and $y$.

Let $\tau_y$ be the first hitting time of $y$ by the random walk:

$$\tau_y \stackrel{\text{def.}}{=} \inf\{k \geq 0 : Z_k = y\}.$$

When $y$ is never attained, let $\tau_y = \infty$. The *hitting probability* of $y$ starting at $x$ is

$$F(x,y) \stackrel{\text{def.}}{=} \mathbb{P}^x[\tau_y < \infty].$$

Note that $F(x,y)$ is positive since the support of $\mu$ generates $\Gamma$, and that $F$ and $G$ are invariant by left diagonal multiplication. In particular, $G(y,y) = G(e,e)$. A straightforward computation (using the strong Markov property) shows that the functions $F$ and $G$ are proportional:

$$(2.1) \qquad G(x,y) = G(y,y)F(x,y) = G(e,e)F(x,y).$$

The metric we will use is the *Green metric* (or hitting metric, defined in [4]):

$$d_G(x,y) \stackrel{\text{def.}}{=} -\ln F(x,y) = \ln G(e,e) - \ln G(x,y).$$

Throughout the article, we will call (with some abuse of notation) *metric* any nonnegative real function $d(\cdot,\cdot)$ on $\Gamma \times \Gamma$ which satisfies the triangle inequality, vanishes on the diagonal and satisfies

$$(2.2) \qquad d(x,y) = 0 = d(y,x) \implies x = y.$$

LEMMA 2.1 ([4], Lemma 2.1). *The function $d_G(\cdot,\cdot)$ is a left-invariant metric on $\Gamma$.*

PROOF. As $F(x,y)$ is bounded by 1, $d_G(\cdot,\cdot)$ is nonnegative. It is also clear that $F(x,x) = 1$ and therefore $d_G(x,x) = 0$ for any $x \in \Gamma$.

The invariance of $F(\cdot,\cdot)$ by left diagonal multiplication implies the same property for $d_G(\cdot,\cdot)$. Also, note that, since the random walk is transient, we have

$$\forall x \neq y \quad 1 > \mathbb{P}^x[\tau'_x < \infty] \geq \mathbb{P}^x[\tau_y < \infty]\mathbb{P}^y[\tau_x < \infty] = F(x,y)F(y,x),$$

where $\tau'_x \stackrel{\text{def.}}{=} \inf\{k \geq 1 : Z_k = x\}$. Thus,

$$d_G(x,y) = d_G(y,x) = 0 \iff F(x,y) = F(y,x) = 1 \iff x = y.$$

Finally,

$$\mathbb{P}^x[\tau_z < \infty] \geq \mathbb{P}^x[\tau_y < \infty]\mathbb{P}^y[\tau_z < \infty]$$

leads to the triangle inequality $d_G(x,z) \leq d_G(x,y) + d_G(y,z)$. □

REMARK 2.2. For the Green metric, and if we assume that $\Gamma$ is not isomorphic to $\mathbb{Z}$, a stronger property than condition (2.2) actually holds, namely,

(2.3) $$d_G(x,y) = 0 \quad \Longrightarrow \quad x = y.$$

The proof of this is as follows. Let $\mathcal{S} \stackrel{\text{def.}}{=} \{x \text{ s.t. } \mu(x) > 0\}$ be the support of the measure $\mu$ and define $\mathcal{E}_0 \stackrel{\text{def.}}{=} \{x \neq e \text{ s.t. } d_G(e,x) = 0\}$. We recall that, as in the Introduction, $\mathcal{S}$ is assumed to generate the whole group $\Gamma$. First, observe that $\mathcal{E}_0$ cannot contain two different elements of $\mathcal{S}$. Indeed, assume that $x \neq y$ both belong to $\mathcal{E}_0 \cap \mathcal{S}$. The random walk has a positive probability of visiting $x$ before $y$. But, since $F(e,y) = 1$, it implies that $F(x,y) = 1$ and therefore $d_G(x,y) = 0$. Likewise, we have $d_G(y,x) = 0$, which contradicts property (2.2). Similar arguments show that if $\mathcal{E}_0$ is not empty, then there exists $a \in \Gamma$ such that $\mathcal{E}_0 = \{a^n; n \in \mathbb{N}^*\}$. One then argues that any element of $\mathcal{S}$ must also be a power of $a$ and, since $\mathcal{S}$ generates $\Gamma$, we conclude that $\Gamma$ is, in fact, the group generated by $a$.

Observe that if $\mu$ is symmetric [$\mu(x) = \mu(x^{-1})$ for all $x \in \Gamma$], then the Green function $G(\cdot,\cdot)$ and the Green metric $d_G$ are also symmetric and therefore $d_G$ becomes a genuine distance on $\Gamma$.

2.2. *Entropy and Green speed.* The measure $\mu$ is now supposed to have finite *entropy*:

$$H(\mu) \stackrel{\text{def.}}{=} -\sum_{x \in \Gamma} \mu(x) \ln \mu(x) < \infty.$$

The first moment of $\mu$ in the Green metric is, by definition, the expected Green distance between $e$ and $Z_1$, which is also the expected Green distance between $Z_n$ and $Z_{n+1}$ for any $n$, having the following analytic expression:

$$\mathbb{E}[d_G(e, Z_1)] = \sum_{x \in \Gamma} \mu(x) \cdot d_G(e, x).$$

LEMMA 2.3. *The finiteness of the entropy $H(\mu)$ implies the finiteness of the first moment of $\mu$ with respect to the Green metric.*

PROOF. By construction, the law of $Z_1 = X_1$ under $\mathbb{P}$ is $\mu$. Since $\mathbb{P}[\tau_x < \infty] \geq \mathbb{P}[Z_1 = x] = \mu(x)$ holds, we have

$$\sum_{x \in \Gamma} \mu(x) \cdot d_G(e, x) = -\sum_{x \in \Gamma} \mu(x) \cdot \ln(\mathbb{P}[\tau_x < \infty])$$
$$\leq -\sum_{x \in \Gamma} \mu(x) \cdot \ln(\mu(x)) = H(\mu)$$

so that $H(\mu) < \infty \Longrightarrow \mathbb{E}[d_G(e, X_1)] < \infty$. □

Let $\ell_G$ be the rate of escape of the random walk $Z_n$ in the Green metric $d_G(e, .)$.

$$\ell_G = \ell_G(\mu) \stackrel{\text{def.}}{=} \lim_{n \to \infty} \frac{d_G(e, Z_n)}{n}$$
$$= \lim_{n \to \infty} \frac{-\ln F(e, Z_n)}{n} = \lim_{n \to \infty} \frac{-\ln G(e, Z_n)}{n},$$

since the functions $F(\cdot, \cdot)$ and $G(\cdot, \cdot)$ differ only by a multiplicative constant. We call $\ell_G$ the *Green speed*. Under the hypothesis that $\mu$ has finite entropy, by the subadditive ergodic theorem (Kingman [22], Derriennic [9]), this limit exists almost surely and in $L^1$.

The sub-additive ergodic theorem of Kingman also allows one to define the asymptotic entropy as the almost sure and $L^1$ limit:

$$h \stackrel{\text{def.}}{=} \lim_{n \to \infty} \frac{-\ln \mu^n(Z_n)}{n},$$

where $\mu^n$ is the $n$th convolution power of the measure $\mu$.

Taking expectations, we deduce that $h$ also satisfies

$$h = \lim_n \frac{H(\mu^n)}{n}.$$

The properties of the asymptotic entropy are studied in great generality in the articles mentioned in the Introduction. In particular, it turns out that $h$ can also be interpreted as a Fisher information. We shall use this fact to conclude the proof of our theorem; see Lemma 2.6.

2.3. *Martin boundary and proof of Theorem* 1.1. The *Martin kernel* is defined [using (2.1)] for all $(x, y) \in \Gamma \times \Gamma$ by

$$K(x, y) \stackrel{\text{def.}}{=} \frac{G(x, y)}{G(e, y)} = \frac{F(x, y)}{F(e, y)}.$$

The Martin kernel continuously extends to a compactification of $\Gamma$ called the *Martin compactification* $\Gamma \cup \partial_M \Gamma$, where $\partial_M \Gamma$ is the *Martin boundary*. Let us briefly recall the construction of $\partial_M \Gamma$. Let $\Psi : \Gamma \to C(\Gamma)$ be defined by $y \longmapsto K(\cdot, y)$. Here, $C(\Gamma)$ is the space of real-valued functions defined on $\Gamma$, endowed with the topology of pointwise convergence. It turns out that $\Psi$ is injective and we may thus identify $\Gamma$ with its image. The closure of $\Psi(\Gamma)$ is compact in $C(\Gamma)$ and, by definition, $\partial_M \Gamma = \overline{\Psi(\Gamma)} \setminus \Psi(\Gamma)$ is the Martin boundary. In the compact space $\Gamma \cup \partial_M \Gamma$, for any initial point $x$, the random walk $Z_n$ almost surely converges to some random variable $Z_\infty \in \partial_M \Gamma$ (see, e.g., Dynkin [12] or Woess [25]).

We note that, by means of the Green metric, one can also consider the Martin compactification as a special example of a *Busemann compactification*. We recall that the Busemann compactification of a proper metric space $(X, d)$ is obtained through the embedding $\Phi \colon X \to C(X)$, defined by $y \longmapsto d(\cdot, y) - d(e, y)$. (Here, $e$ denotes an arbitrary base point.) In general, $C(X)$ should be endowed with the topology of uniform convergence on compact sets. The Busemann compactification of $X$ is the closure of the image $\Phi(X)$ in $C(X)$. We refer to Ballmann, Gromov and Schroeder [2] and to Karlsson and Ladrappier [20] and the references therein for further details.

If one now chooses for $X$ the group $\Gamma$ itself and, for the distance $d$, the Green metric, constructions of both the Martin and Busemann compactifications coincide, as is straightforward from the relation

$$d_G(\cdot, y) - d_G(e, y) = -\ln K(\cdot, y).$$

We first prove that the Green speed can be expressed in terms of the extended Martin kernel. Theorem 1.1 will then be a direct consequence of the formulas in Proposition 2.4 and Lemma 2.6. For that purpose, we need to define the reversed law $\tilde{\mu}$:

$$\forall x \in \Gamma \qquad \tilde{\mu}(x) \stackrel{\text{def.}}{=} \mu(x^{-1}).$$

Note that $H(\tilde{\mu}) = H(\mu)$.

PROPOSITION 2.4. *Let $\mu$ be a probability measure on $\Gamma$ with finite entropy $H(\mu)$ and whose support generates $\Gamma$. Let $(Z_n)$ be a random walk on $\Gamma$ of law $\mu$ (starting at $e$) and let $\tilde{X}_1$ be an independent random variable of law $\tilde{\mu}$. Then,*

$$\ell_G = \mathbb{E}\tilde{\mathbb{E}}[-\ln K(\tilde{X}_1, Z_\infty)],$$

*where $\tilde{\mathbb{E}}$ refers to the integration with respect to the random variable $\tilde{X}_1$ and $\mathbb{E}$ refers to the integration with respect to the random walk $(Z_n)$.*

PROOF. As $\mu$ is supposed to have finite entropy, $\ell_G$ is well defined as an almost sure and $L^1$ limit. We will prove that the sequence

$$\mathbb{E}[d_G(e, Z_{n+1}) - d_G(e, Z_n)] = \mathbb{E}[-\ln G(e, Z_{n+1}) + \ln G(e, Z_n)]$$

converges to $\mathbb{E}\tilde{\mathbb{E}}[-\ln K(\tilde{X}_1, Z_\infty)]$. Since its limit in the Cesaro sense is $\ell_G$, it implies the formula in Proposition 2.4.

By definition of the reversed law $\tilde{\mu}$, $\tilde{X}_1^{-1}$ has the same law as $X_1$, the first increment of the random walk $(Z_n)$. Also, note that $X_2 \cdots X_{n+1}$ has the same law as $Z_n = X_1 \cdots X_n$. Since we have assumed that $\tilde{X}_1$ is independent of the sequence $(Z_n)$, $Z_{n+1} = X_1 \cdot X_2 \cdots X_{n+1}$ has the same law as $\tilde{X}_1^{-1} \cdot Z_n$ and therefore, using the translation invariance, $G(e, Z_{n+1})$ has the same law as $G(\tilde{X}_1, Z_n)$. Thus,

$$\mathbb{E}[-\ln G(e, Z_{n+1}) + \ln G(e, Z_n)] = \mathbb{E}\tilde{\mathbb{E}}[-\ln G(\tilde{X}_1, Z_n) + \ln G(e, Z_n)]$$
$$= \mathbb{E}\tilde{\mathbb{E}}[-\ln K(\tilde{X}_1, Z_n)].$$

By continuity of the Martin kernel up to the Martin boundary, for every $x \in \Gamma$, the sequence $K(x, Z_n)$ almost surely converges to $K(x, Z_\infty)$. We need an integrable bound for $-\ln K(\tilde{X}_1, Z_n)$ (uniformly in $n$) to justify the convergence of the expectation.

To prove that $-\ln K(\tilde{X}_1, Z_n)$ cannot go too far in the negative direction, we first prove a maximal inequality for the sequence $(K(\tilde{X}_1, Z_n))_n$, following Dynkin [12].

LEMMA 2.5. *For any $a > 0$,*

$$\mathbb{P}\tilde{\mathbb{P}}\left[\sup_n K(\tilde{X}_1, Z_n) \geq a\right] \leq \frac{1}{a},$$

*where $\tilde{\mathbb{P}}$ refers with the measure associated with the random variable $\tilde{X}_1$ and $\mathbb{P}$ refers to the measure associated with the random walk $(Z_n)$.*

PROOF. We fix an integer $R$. Let $\sigma_R$ be the time of the last visit to the ball $B_G(e, R) \stackrel{\text{def.}}{=} \{x \in \Gamma \text{ s.t. } d_G(e, x) \leq R\}$ for the random walk $(Z_n)$. [Since the random walk is transient, $\sigma_R$ is well defined and almost surely finite if the random walk starts within $B_G(e, R)$. Otherwise, $\sigma_R$ is set to be infinite when $B_G(e, R)$ is never reached.] Let us define the sequence $(Z_{\sigma_R - k})$ ($k \in \mathbb{N}$). As this sequence (in $\Gamma$) is only defined for $k \leq \sigma_R$, we take the following convention for negative indices:

$$\{k > \sigma_R\} \implies \{Z_{\sigma_R - k} \stackrel{\text{def.}}{=} \star\}.$$

In this way, the sequence $(Z_{\sigma_R - k})_{k \in \mathbb{N}}$ is well defined and takes its values in $\Gamma \cup \{\star\}$. Note that $Z_{\sigma_R}$ takes its value in $B_G(e, R)$.

Let us call $\mathcal{F}_k$ the $\sigma$-algebra generated by $(Z_{\sigma_R}, \ldots, Z_{\sigma_R - k})$ and observe that

$$\mathbb{1}_{\{k \leq \sigma_R\}} \in \mathcal{F}_k$$

since $\{k \leq \sigma_R\}$ means that none of $Z_{\sigma_R}, \ldots, Z_{\sigma_R-k}$ equals $\star$. With the convention that, for any $x \in \Gamma$, $K(x, \star) = 0$, we can define, for any $x$ in $\Gamma$, the nonnegative sequence $(K(x, Z_{\sigma_R-k}))$ $(k \in \mathbb{N})$. This sequence is adapted to the filtration $(\mathcal{F}_k)$ and we will prove, following Dynkin [12], Sections 6, 7 that it is a supermartingale with respect to $(\mathcal{F}_k)$.

For this purpose, let us check that, for any positive integer $k$ and any sequence $z_0, z_1, \ldots, z_{k-1}$ in $\Gamma \cup \{\star\}$ [with $z_0 \in B_G(e, R)$],

$$\mathbb{E}\left[K(x, Z_{\sigma_R-k}) \prod_{j=0}^{k-1} \mathbb{1}_{\{Z_{\sigma_R-j}=z_j\}}\right]$$
(2.4)
$$= (K(x, z_{k-1}) - \delta_x(z_{k-1})G(e, x)^{-1}) \cdot \mathbb{E}\left[\prod_{j=0}^{k-1} \mathbb{1}_{\{Z_{\sigma_R-j}=z_j\}}\right].$$

We first compute the left-hand side of (2.4) in the case where none of $z_0, z_1, \ldots, z_{k-1}$ equals $\star$. First using the fact that $K(x, \star) = 0$, we have

$$\sum_{z_k \in \Gamma \cup \{\star\}} \mathbb{P}[Z_{\sigma_R} = z_0, \ldots, Z_{\sigma_R-(k-1)} = z_{k-1}, Z_{\sigma_R-k} = z_k] \cdot K(x, z_k)$$

$$= \sum_{z_k \in \Gamma} \mathbb{P}[Z_{\sigma_R} = z_0, \ldots, Z_{\sigma_R-k} = z_k] \cdot K(x, z_k)$$

$$= \sum_{z_k \in \Gamma} \mathbb{P}[k \leq \sigma_R, Z_{\sigma_R} = z_0, \ldots, Z_{\sigma_R-k} = z_k] \cdot K(x, z_k)$$

since the fact that none of $z_0, \ldots, z_k$ equals $\star$ means, in particular, that

$$\bigcap_{j=0}^{k} \{Z_{\sigma_R-j} = z_j\} \subset \{k \leq \sigma_R\}.$$

Then,

$$\sum_{z_k \in \Gamma \cup \{\star\}} \mathbb{P}[Z_{\sigma_R} = z_0, \ldots, Z_{\sigma_R-(k-1)} = z_{k-1}, Z_{\sigma_R-k} = z_k] \cdot K(x, z_k)$$

$$= \sum_{z_k \in \Gamma} \sum_{m=k}^{\infty} \mathbb{P}[\sigma_R = m, Z_m = z_0, \ldots, Z_{m-k} = z_k] \cdot K(x, z_k)$$

$$= \sum_{z_k \in \Gamma} \sum_{m=k}^{\infty} \mathbb{P}[Z_{m-k} = z_k]\mu(z_k^{-1}z_{k-1})\cdots\mu(z_1^{-1}z_0)\mathbb{P}^{z_0}[\sigma_R = 0] \cdot K(x, z_k)$$

$$= \mu(z_{k-1}^{-1}z_{k-2})\cdots\mu(z_1^{-1}z_0)\mathbb{P}^{z_0}[\sigma_R = 0] \sum_{z_k \in \Gamma} G(e, z_k)\mu(z_k^{-1}z_{k-1}) \cdot K(x, z_k)$$

$$= \mu(z_{k-1}^{-1}z_{k-2})\cdots\mu(z_1^{-1}z_0)\mathbb{P}^{z_0}[\sigma_R = 0] \sum_{z_k \in \Gamma} G(x, z_k)\mu(z_k^{-1}z_{k-1})$$

$$= \mu(z_{k-1}^{-1}z_{k-2}) \cdots \mu(z_1^{-1}z_0)\mathbb{P}^{z_0}[\sigma_R = 0](G(x, z_{k-1}) - \delta_x(z_{k-1})).$$

Using the same kind of computation, we get that the right-hand side of (2.4) equals

$$\sum_{m=k-1}^{\infty} \mathbb{P}[\sigma_R = m, Z_m = z_0, \ldots, Z_{m-(k-1)} = z_{k-1}]$$

$$\times (K(x, z_{k-1}) - \delta_x(z_{k-1})G(e, x)^{-1})$$

$$= \sum_{m=k-1}^{\infty} \mathbb{P}[Z_{m-(k-1)} = z_{k-1}]\mu(z_{k-1}^{-1}z_{k-2}) \cdots \mu(z_1^{-1}z_0)\mathbb{P}^{z_0}[\sigma_R = 0]$$

$$\times (K(x, z_{k-1}) - \delta_x(z_{k-1})G(e, x)^{-1})$$

$$= \mu(z_{k-1}^{-1}z_{k-2}) \cdots \mu(z_1^{-1}z_0)\mathbb{P}^{z_0}[\sigma_R = 0](G(x, z_{k-1}) - \delta_x(z_{k-1})).$$

So, (2.4) is true as soon as $z_0, \ldots, z_{k-1}$ take values in $\Gamma$. Now, supposing that $z_j = \star$ for some $j \leq k-1$, we have

$$\{Z_{\sigma_R - j} = z_j\} \implies \{Z_{\sigma_R - (k-1)} = \star\} \implies \{Z_{\sigma_R - k} = \star\}.$$

Since $K(x, \star) = 0$, the left-hand side of (2.4) is zero. To check that the right-hand side is also zero, observe that

$$z_{k-1} \neq \star \implies \mathbb{1}_{\{Z_{\sigma_R - j} = z_j\}} \cdot \mathbb{1}_{\{Z_{\sigma_R - (k-1)} = z_{(k-1)}\}} = 0$$

$$\implies \mathbb{E}\left[\prod_{j=0}^{k-1} \mathbb{1}_{\{Z_{\sigma_R - j} = z_j\}}\right] = 0,$$

and, as $x \in \Gamma$,

$$z_{k-1} = \star \implies K(x, z_{k-1}) = 0 \quad \text{and} \quad \delta_x(z_{k-1}) = 0.$$

The proof of (2.4) is now complete. Since the Green function is positive, we deduce from (2.4) that

$$\mathbb{E}\left[K(x, Z_{\sigma_R - k})\prod_{j=0}^{k-1} \mathbb{1}_{\{Z_{\sigma_R - j} = z_j\}}\right] \leq K(x, z_{k-1}) \cdot \mathbb{E}\left[\prod_{j=0}^{k-1} \mathbb{1}_{\{Z_{\sigma_R - j} = z_j\}}\right],$$

thus proving the supermartingale property of the sequence $(K(x, Z_{\sigma_R - k}))$ ($k \in \mathbb{N}$).

We use similar arguments to compute the expectation of the value of the supermartingale at time $k = 0$, $\mathbb{E}[K(x, Z_{\sigma_R})]$, which turns out not to depend on $R$:

$$\mathbb{E}[K(x, Z_{\sigma_R})] = \sum_{m=0}^{\infty} \sum_{z \in B_G(e, R)} \mathbb{P}[\sigma_R = m, Z_m = z] \cdot K(x, z)$$

$$= \sum_{m=0}^{\infty} \sum_{z \in B_G(e,R)} \mathbb{P}^z[\sigma_R = 0] \cdot \mathbb{P}[Z_m = z] \cdot K(x,z)$$

$$= \sum_{z \in B_G(e,R)} \mathbb{P}^z[\sigma_R = 0] \cdot G(x,z)$$

$$= \sum_{z \in B_G(e,R)} \mathbb{P}^z[\sigma_R = 0] \sum_{m=0}^{\infty} \mathbb{P}^x[Z_m = z]$$

$$= \sum_{z \in B_G(e,R)} \sum_{m=0}^{\infty} \mathbb{P}^x[\sigma_R = m, Z_{\sigma_R} = z]$$

$$= \mathbb{P}^x[\sigma_R < \infty] \leq 1.$$

We can now use Doob's maximal inequality for nonnegative supermartingales (see, e.g., Breiman [6], Proposition 5.13) to get that

$$\forall x \in \Gamma \quad \mathbb{P}\left[\sup_k K(x, Z_{\sigma_R - k}) \geq a\right] \leq \frac{1}{a}.$$

So, $\mathbb{P}\tilde{\mathbb{P}}[\sup_k K(\tilde{X}_1, Z_{\sigma_R - k}) \geq a] \leq \frac{1}{a}$ and, letting $R$ tend to infinity,

$$\mathbb{P}\tilde{\mathbb{P}}\left[\sup_n K(\tilde{X}_1, Z_n) \geq a\right] \leq \frac{1}{a}. \qquad \Box$$

Let us return to the proof of Proposition 2.4. Lemma 2.5 implies that, for any $b > 0$,

$$\mathbb{P}\tilde{\mathbb{P}}\left[\sup_n \ln K(\tilde{X}_1, Z_n) \geq b\right] \leq e^{-b}$$

and therefore $\mathbb{E}\tilde{\mathbb{E}}[\sup_n \ln K(\tilde{X}_1, Z_n) \mathbb{1}_{K(\tilde{X}_1, Z_n) \geq 1}] < \infty$.

On the other hand, we have

$$K(x, Z_n) = \frac{\mathbb{P}^x[\tau_{Z_n} < \infty]}{\mathbb{P}^e[\tau_{Z_n} < \infty]} \geq \frac{\mathbb{P}^x[\tau_e < \infty] \cdot \mathbb{P}^e[\tau_{Z_n} < \infty]}{\mathbb{P}^e[\tau_{Z_n} < \infty]}$$

$$= \mathbb{P}^e[\tau_{x^{-1}} < \infty] \geq \tilde{\mu}(x)$$

and

$$\tilde{\mathbb{E}}[-\ln \tilde{\mu}(\tilde{X}_1)] = H(\tilde{\mu}) = H(\mu) < \infty.$$

Writing that

$$|\ln K(\tilde{X}_1, Z_n)| = \ln K(\tilde{X}_1, Z_n) \mathbb{1}_{K(\tilde{X}_1, Z_n) \geq 1} - \ln K(\tilde{X}_1, Z_n) \mathbb{1}_{K(\tilde{X}_1, Z_n) \leq 1}$$

$$\leq \ln K(\tilde{X}_1, Z_n) \mathbb{1}_{K(\tilde{X}_1, Z_n) \geq 1} - \ln \tilde{\mu}(\tilde{X}_1),$$

we conclude that the random variable $\sup_n |\ln K(\tilde{X}_1, Z_n)|$ is integrable. We can therefore apply the dominated convergence theorem to deduce that the sequence $\mathbb{E}[-\ln G(e, Z_{n+1}) + \ln G(e, Z_n)]$ converges to

$$\mathbb{E}\tilde{\mathbb{E}}[-\ln K(\tilde{X}_1, Z_\infty)]. \qquad \square$$

LEMMA 2.6. *Let $\Gamma$ be a countable group and $\mu$ be a probability measure on $\Gamma$ whose support generates $\Gamma$ with finite entropy $H(\mu)$. Then,*

$$h = \mathbb{E}\tilde{\mathbb{E}}[-\ln K(\tilde{X}_1, Z_\infty)].$$

PROOF. Recall that $\tilde{\mu}$ is the law of $\tilde{X}_1$. We have

$$\mathbb{E}\tilde{\mathbb{E}}[-\ln K(\tilde{X}_1, Z_\infty)] = \int_\Gamma \int_{\partial_M \Gamma} -\ln(K(x, \xi)) \, d\nu(\xi) \, d\tilde{\mu}(x),$$

where $\nu_y(\cdot)$ is the harmonic measure on the Martin boundary $\partial_M \Gamma$ for a random walk (of law $\mu$) starting at $y$ and $\nu(\cdot) = \nu_e(\cdot)$. By the Martin boundary convergence theorem (see Hunt [16] or Woess [25], Theorem 24.10) the Martin kernel $K(x, \xi)$ is the Radon–Nikodym derivative of $\nu_x$ by $\nu$ at $\xi$. Therefore,

$$\mathbb{E}\tilde{\mathbb{E}}[-\ln K(\tilde{X}_1, Z_\infty)] = \int_\Gamma \int_{\partial_M \Gamma} -\ln\left(\frac{d\nu_x(\xi)}{d\nu(\xi)}\right) d\nu(\xi) \, d\mu(x^{-1}).$$

We will make the following changes of variables. As $\partial_M \Gamma$ is stable by left multiplication, the change of variables $\xi \longmapsto x^{-1}\xi$ gives $\nu_x(\xi) \longmapsto \nu(\xi)$ and $\nu(\xi) \longmapsto \nu_{x^{-1}}(\xi)$. Hence, also changing $x$ into $x^{-1}$ gives

$$\mathbb{E}\tilde{\mathbb{E}}[-\ln K(\tilde{X}_1, Z_\infty)] = \int_\Gamma \int_{\partial_M \Gamma} -\ln\left(\frac{d\nu(\xi)}{d\nu_x(\xi)}\right) d\nu_x(\xi) \, d\mu(x)$$

(2.5)

$$= \int_\Gamma \int_{\partial_M \Gamma} \ln\left(\frac{d\nu_x(\xi)}{d\nu(\xi)}\right) d\nu_x(\xi) \, d\mu(x).$$

Observe that $d\nu_x(\xi)/d\nu(\xi)$ is the Radon–Nikodym derivative of the joint law of $(\tilde{X}_1^{-1}, Z_\infty)$ with respect to the product measure $\mu(\cdot) \otimes \nu(\cdot)$. Therefore, (2.5) means that $\mathbb{E}\tilde{\mathbb{E}}[-\ln K(\tilde{X}_1, Z_\infty)]$ is the relative entropy of the joint law of $(\tilde{X}_1^{-1}, Z_\infty)$ with respect to $\mu(\cdot) \otimes \nu(\cdot)$, which equals the asymptotic entropy $h$ (see Derriennic [11], who actually takes the latter as the definition of the asymptotic entropy and proves that the two definitions coincide). $\square$

**3. Finitely generated groups.** We now restrict ourselves to a finitely generated group $\Gamma$.

3.1. *Volume growth in the Green metric.* For a given finite generating set $\mathcal{S}$, we define the associated *word metric*:

$$d_w(x,y) \stackrel{\text{def.}}{=} \min\{n \text{ s.t. } x^{-1}y = g_1 g_2 \cdots g_n \text{ with } g_i \in \mathcal{S}\}.$$

This distance is the geodesic graph distance of the Cayley graph of $\Gamma$ defined by $\mathcal{S}$. Different choices of generating sets lead to different word distances in the same quasi-isometry class. When $\mu$ is symmetric and finitely supported, the two metrics $d_G(\cdot,\cdot)$ and $d_w(\cdot,\cdot)$ can be compared (see [4], Lemma 2.2). These two metrics are equivalent for any nonamenable group and also for some amenable groups (e.g., the Lamplighter group $\mathbb{Z} \wr \mathbb{Z}_2$).

Throughout the article, the notion of growth of the group $\Gamma$ always refers to the function $V_w(n) \stackrel{\text{def.}}{=} \#\{x \in \Gamma \text{ s.t. } d_w(e,x) \leq n\}$ for some (equivalently, any) symmetric finite generating set. The group will be said to have:

- polynomial growth when $V_w(n) = O(n^D)$ for some constant $D$ (the largest integer $D$ satisfying this condition is called the *degree* of the group);
- superpolynomial growth when $V_w(n)/n^D$ tends to infinity for every $D$;
- subexponential growth when $V_w(n) = o(e^{Cn})$ for every constant $C > 0$;
- exponential growth when $V_w(n)/e^{Cn}$ tends to infinity for some $C > 0$.

We are now interested in the asymptotic behavior of the volume of the balls for the Green metric. Let us define $B_G(e,n) \stackrel{\text{def.}}{=} \{x \in \Gamma \text{ s.t. } d_G(e,x) \leq n\}$, $V_G(n) \stackrel{\text{def.}}{=} \#B_G(e,n)$ and the corresponding logarithmic volume growth,

$$v_G \stackrel{\text{def.}}{=} \limsup_{n \to \infty} \frac{\ln(V_G(n))}{n}.$$

PROPOSITION 3.1. *Let us suppose that $\Gamma$ is not a finite extension of $\mathbb{Z}$ or $\mathbb{Z}^2$. For any random walk on $\Gamma$:*

  (i) *if $\Gamma$ has superpolynomial growth, then $v_G \leq 1$;*
  (ii) *if $\Gamma$ has polynomial growth of degree $D$, then $v_G \leq \frac{D}{D-2}$.*

PROOF. Observe that Proposition 2.3 in [4] proves (i) when $\mu$ has finite support and is symmetric.

We recall the following classical result (see, e.g., Woess [25]). Let $\mu$ be a symmetric measure with finite support and let $\Gamma$ have at least polynomial growth of degree $D$ ($D \geq 3$). Then,

(3.1) $\quad \exists C_e > 1 \text{ s.t. } \forall x, y \in \Gamma \text{ and } k \in \mathbb{N} \qquad \mathbb{P}^x[Z_k = y] \leq C_e k^{-D/2}.$

The above estimate remains valid even without the symmetry and the finite support hypotheses. Indeed, Coulhon's result ([7], Proposition IV.4) (see also Coulhon and Saloff-Coste [8]) allows one to extend upper bounds of the $n$th convolution power of a symmetric probability measure $\mu_1$ to the $n$th convolution power of another probability measure $\mu_2$ under the following condition:

(3.2) $$\exists c > 0 \text{ s.t. } \forall x \quad \mu_1(x) \leq c\mu_2(x).$$

For a general probability measure $\mu$ whose support generates $\Gamma$, there exists $K$ such that the support of $\mu^K$ contains any finite symmetric generating set $\mathcal{S}$ of $\Gamma$.

Hence, choosing $\mu_2 = \mu^K$, $c = (\min_{x \in \mathcal{S}} \mu_2(x))^{-1}$ and $\mu_1 = (1/\#\mathcal{S}) \times \delta_\mathcal{S}(x)$, the uniform distribution on $\mathcal{S}$, we see that the measures $\mu_1$ and $\mu_2$ satisfy condition (3.2). Therefore, the estimate (3.1) remains valid for $\mu$, with a possibly different constant $C_e$.

The same argument as in [4] shows that (3.1) implies

$$V_G(n) \leq C \exp\left(\frac{D}{D-2} \cdot n\right)$$

for some constant $C$. Thus, $v_G \leq \frac{D}{D-2}$. For groups with superpolynomial growth, letting $D$ go to infinity gives $v_G \leq 1$. □

REMARK 3.2. If the measure $\mu$ has finite support, then it is already known that $v_G \geq 1$ [4] Proposition 2.3. From Lemma 3.3 and Proposition 3.4, we will also get that $v_G \geq 1$ when $\mu$ has finite entropy and $h > 0$, but $\mu$ may have an infinite support. This implies that $v_G = 1$ for groups with superpolynomial growth and measures of finite entropy such that $h > 0$.

3.2. *The "fundamental" inequality.* We now present a different proof of Theorem 1.1 in the case of finitely generated groups. The interest of this proof comes from an extended version of the "fundamental" inequality relating the asymptotic entropy, the logarithmic volume growth and the rate of escape.

There is the following general, obvious link between the Green speed and the asymptotic entropy.

LEMMA 3.3. *For any random walk with finite entropy $H(\mu)$, we have $\ell_G \leq h$.*

PROOF. The sequence $\frac{1}{n} d_G(e, Z_n)$ converges to $\ell_G$ in $L^1$. Therefore,

$$\ell_G = \lim_{n \to \infty} \frac{-\sum_{x \in \Gamma} \mu^n(x) \ln(\sum_{k=0}^{\infty} \mu^k(x))}{n}$$
$$\leq \lim_{n \to \infty} \frac{-\sum_{x \in \Gamma} \mu^n(x) \ln \mu^n(x)}{n} = h. \qquad \Box$$

Our aim is to prove the other inequality and deduce that $h = \ell_G$.

*Groups with polynomial volume growth.* For groups with polynomial growth, Lemma 3.3 gives the (trivial) equality since any random walk has zero asymptotic entropy. Indeed, these groups have a trivial Poisson boundary (Dynkin and Malyutov [13]), which is equivalent to $h = 0$ for measures with finite entropy (Derriennic [10] and Kaimanovich and Vershik [18], see also Kaimanovich [17], Theorem 1.6.7).

*Groups with superpolynomial volume growth.* We rely on the so-called fundamental inequality.

$$h \leq \ell_G \cdot v_G, \tag{3.3}$$

which holds when $\mu$ has finite entropy. For groups with superpolynomial growth, Proposition 3.1 gives $v_G \leq 1$, therefore inequality (3.3) implies that $h \leq \ell_G$ and we conclude that $h = \ell_G$. Thus, all that remains to be done in order to complete the proof of Theorem 1.1 in the case of groups with superpolynomial growth is to justify (3.3). This is the content of the next proposition.

A version of inequality (3.3), when the speed and volume growth are computed in a word metric, is proved by Guivarc'h [15] and discussed in great detail by Vershik [24]. The same proofs as in [15] or [24] would apply to any invariant metric on $\Gamma$, for instance, the Green metric, the provided $\mu$ has finite support. The fundamental inequality is also known to hold for measures with unbounded support and a finite first moment in a word metric. See, for instance, Erschler [14], Lemma 6 or Karlsson and Ledrappier [21], but note that their argument seems to apply only to word metrics and observe that the Green metric is not a word metric in general (as a matter of fact, it need not even be a geodesic metric). We shall derive the fundamental inequality in the Green metric, under the simple assumption that the entropy of $\mu$ is finite.

We present our result in a general setting (for any invariant metric and group) since it has its own interest.

PROPOSITION 3.4. *Let $\mu$ be the law of the increment of a random walk on a countable group $\Gamma$, starting at a point $e$, and let $d(\cdot, \cdot)$ be a left-invariant metric. Under the hypothesis that:*

- *the measure $\mu$ has finite entropy;*
- *the measure $\mu$ has finite first moment with respect to the metric $d$;*
- *the logarithmic volume growth $v \stackrel{\text{def.}}{=} \limsup_{n \to \infty} \frac{\ln(\#B(e,n))}{n}$ is finite, then*

*the asymptotic entropy $h$, the rate of escape $\ell \stackrel{\text{def.}}{=} \lim_n \frac{d(e,Z_n)}{n}$ (limit both in $L^1$ and almost surely) and the logarithmic volume growth $v$ satisfy the following inequality:*

$$h \leq \ell \cdot v.$$

PROOF. The proof relies on an idea of Guivarc'h [15] Proposition C.2. Fix $\varepsilon > 0$ and, for all integers $n$, let $B_\varepsilon^n \stackrel{\text{def.}}{=} B(e, (\ell+\varepsilon)n)$ [here, the balls are defined for the metric $d(e, \cdot)$]. We split $\Gamma \setminus B_\varepsilon^n$ into a sequence of annuli as follows. Choose $K > \ell + \varepsilon$ and define

$$\mathcal{C}_\varepsilon^{n,K} \stackrel{\text{def.}}{=} B(e, Kn) \setminus B_\varepsilon^n,$$

$$\forall i \geq 1 \quad \mathcal{C}_i^{n,K} \stackrel{\text{def.}}{=} B(e, 2^i Kn) \setminus B(e, 2^{i-1} Kn).$$

Define the conditional entropy

$$H(\mu|A) \stackrel{\text{def.}}{=} -\sum_{x \in A} \frac{\mu(x)}{\mu(A)} \ln \frac{\mu(x)}{\mu(A)}.$$

The entropy of $\mu^n$ can then be written as

$$\begin{aligned}
H(\mu^n) &= \mu^n(B_\varepsilon^n) \cdot H(\mu^n | B_\varepsilon^n) \\
&\quad + \mu^n(\mathcal{C}_\varepsilon^{n,K}) \cdot H(\mu^n | \mathcal{C}_\varepsilon^{n,K}) \\
&\quad + \sum_{i=1}^\infty \mu^n(\mathcal{C}_i^{n,K}) \cdot H(\mu^n | \mathcal{C}_i^{n,K}) + H_n',
\end{aligned} \tag{3.4}$$

where

$$\begin{aligned}
H_n' &\stackrel{\text{def.}}{=} -\mu^n(B_\varepsilon^n) \cdot \ln(\mu^n(B_\varepsilon^n)) \\
&\quad - \mu^n(\mathcal{C}_\varepsilon^{n,K}) \cdot \ln(\mu^n(\mathcal{C}_\varepsilon^{n,K})) \\
&\quad - \sum_{i=1}^\infty \mu^n(\mathcal{C}_i^{n,K}) \cdot \ln(\mu^n(\mathcal{C}_i^{n,K})).
\end{aligned} \tag{3.5}$$

We will repeatedly use the fact that the entropy of any probability measure supported by a finite set is maximal for the uniform measure and then equals the logarithm of the volume. First, observe that

$$H(\mu^n | B_\varepsilon^n) \leq \ln(\# B_\varepsilon^n) \leq (\ell + \varepsilon) \cdot v \cdot n + o(n)$$

and thus the first term in (3.4) satisfies

$$\lim_n \frac{\mu^n(B_\varepsilon^n) \cdot H(\mu^n | B_\varepsilon^n)}{n} \leq (\ell + \varepsilon) \cdot v.$$

For the second term in (3.4), we get that

$$H(\mu^n | \mathcal{C}_\varepsilon^{n,K}) \leq \ln(\# \mathcal{C}_\varepsilon^{n,K}) \leq K \cdot v \cdot n + o(n).$$

On the other hand, $\ell$ is also the limit in probability of $d(e, Z_n)/n$, hence $\forall \varepsilon > 0$, $\lim_n \mu^n(B_\varepsilon^n) = 1$. Therefore, $\lim_n \mu^n(\mathcal{C}_\varepsilon^{n,K}) = 0$ and the second term in (3.4) satisfies

$$\lim_n \frac{\mu^n(\mathcal{C}_\varepsilon^{n,K}) \cdot H(\mu^n | \mathcal{C}_\varepsilon^{n,K})}{n} = 0.$$

For the third term in (3.4), as before, we have
$$H(\mu^n|\mathcal{C}_i^{n,K}) \leq \ln(\#\mathcal{C}_i^{n,K}) \leq 2^i K \cdot v \cdot n + o(n)$$
and, by the definition of $\mathcal{C}_i^{n,K}$,

(3.6) $$\mu^n(\mathcal{C}_i^{n,K}) = \mathbb{E}[\mathbb{1}_{\{Z_n \in \mathcal{C}_i^{n,K}\}}] \leq \mathbb{E}\left[\frac{d(e, Z_n)}{2^{i-1}Kn} \cdot \mathbb{1}_{\{Z_n \in \mathcal{C}_i^{n,K}\}}\right].$$

So,
$$\frac{1}{n}\sum_{i=1}^{\infty} \mu^n(\mathcal{C}_i^{n,K}) \cdot H(\mu^n|\mathcal{C}_i^{n,K})$$
$$\leq \left(\frac{2v}{n} + o\left(\frac{1}{n}\right)\right) \mathbb{E}\left[d(e, Z_n) \sum_{i=1}^{\infty} \mathbb{1}_{\{Z_n \in \mathcal{C}_i^{n,K}\}}\right]$$
$$= \left(\frac{2v}{n} + o\left(\frac{1}{n}\right)\right) \mathbb{E}[d(e, Z_n) \cdot \mathbb{1}_{\{d(e, Z_n) > Kn\}}].$$

As $d(e, Z_n) \leq \sum_{k=1}^{n} d(e, X_k)$, we have
$$\frac{1}{n}\sum_{i=1}^{\infty} \mu^n(\mathcal{C}_i^{n,K}) \cdot H(\mu^n|\mathcal{C}_i^{n,K})$$
$$\leq \left(\frac{2v}{n} + o\left(\frac{1}{n}\right)\right) \sum_{j=1}^{n} \mathbb{E}[d(e, X_j) \cdot \mathbb{1}_{\{\sum_{k=1}^n d(e, X_k) > Kn\}}]$$
$$= (2v + o(1))\mathbb{E}[d(e, X_1) \cdot \mathbb{1}_{\{\sum_{k=1}^n d(e, X_k) > Kn\}}]$$

since $X_1, \ldots, X_n$ are i.i.d., so the random variables
$$Y_j \stackrel{\text{def.}}{=} d(e, X_j) \cdot \mathbb{1}_{\{\sum_{k=1}^n d(e, X_k) > Kn\}}$$
have the same distribution.

By the strong law of large numbers, the sequence $\frac{1}{n}\sum_{k=1}^n d(e, X_k)$ almost surely converges to $\mathbb{E}[d(e, X_1)] \stackrel{\text{def.}}{=} m < \infty$. As a consequence, for any $K > m$, we have

(3.7) $$d(e, X_1) \cdot \mathbb{1}_{\{\sum_{k=1}^n d(e, X_k) > Kn\}} \xrightarrow{\text{a.s.}} 0.$$

Moreover, as
$$d(e, X_1) \cdot \mathbb{1}_{\{\sum_{k=1}^n d(e, X_k) > Kn\}} \leq d(e, X_1),$$
which is integrable, the limit in (3.7) also occurs in $L^1$. Then,
$$\lim_n \frac{1}{n}\sum_{i=1}^{\infty} \mu^n(\mathcal{C}_i^{n,K}) \cdot H(\mu^n|\mathcal{C}_i^{n,K}) = 0.$$

We are left with $H'_n$. As $\lim_n \mu^n(B_\varepsilon^n) = 1$ and $\lim_n \mu^n(\mathcal{C}_\varepsilon^{n,K}) = 0$, we have

$$\lim_n [-\mu^n(B_\varepsilon^n) \cdot \ln(\mu^n(B_\varepsilon^n)) - \mu^n(\mathcal{C}_\varepsilon^{n,K}) \cdot \ln(\mu^n(\mathcal{C}_\varepsilon^{n,K}))] = 0.$$

For the last term in (3.5), note that (3.6) gives

$$\mu^n(\mathcal{C}_i^{n,K}) \leq \frac{1}{2^{i-1}Kn} \sum_{k=1}^n \mathbb{E}[d(e, X_k)] \leq \frac{m}{2^{i-1}K}.$$

Together with the inequality $-a\ln(a) \leq 2e^{-1}\sqrt{a}$, we get

$$-\sum_{i=1}^\infty \mu^n(\mathcal{C}_i^{n,K}) \cdot \ln(\mu^n(\mathcal{C}_i^{n,K})) \leq 2e^{-1} \sum_{i=1}^\infty \sqrt{\mu^n(\mathcal{C}_i^{n,K})} < \infty.$$

So, $\lim_n H'_n/n = 0$.

Finally, taking the limit $n \to \infty$, we deduce from (3.4) that $h \leq (\ell + \varepsilon) \cdot v$ for any $\varepsilon$, so $h \leq \ell \cdot v$. □

We conclude with a final remark.

REMARK 3.5. The proof of Theorem 1.1 using the Martin boundary relies on the translation invariance of $\Gamma$, but the hypothesis that the graph is a Cayley graph of a countable group seems too strong. It would be interesting to extend this proof to the case of space homogeneous Markov chains (see Kaimanovich and Woess [19]).

**Acknowledgments.** The authors would like to thank the participants of the working group "Boundaries of groups" in Marseille, where fruitful discussions took place. We are also grateful to Yuval Peres for pointing out the reference [3] after a first version of the present article was made public.

S. BLACHÈRE  
EURANDOM  
P.O. BOX 513  
5600 MB EINDHOVEN  
THE NETHERLANDS  
E-MAIL: blachere@eurandom.tue.nl  
URL: http://euridice.tue.nl/~blachere/

P. HAÏSSINSKY  
P. MATHIEU  
CENTRE DE MATHÉMATIQUES  
ET D'INFORMATIQUE UNIVERSITÉ AIX-MARSEILLE 1  
39 RUE JOLIOT-CURIE  
13453 MARSEILLE CEDEX 13  
FRANCE  
E-MAIL: phaissin@cmi.univ-mrs.fr  
pmathieu@cmi.univ-mrs.fr  
URL: http://www.cmi.univ-mrs.fr/~pmathieu/